\documentclass[a4paper,12pt]{article}
\usepackage{amsthm,amsmath}
\usepackage{amssymb}
\usepackage[latin1]{inputenc}
\usepackage[all]{xy}
\usepackage{enumerate}
\usepackage[T1]{fontenc}
\usepackage{a4wide}
\usepackage[mathscr]{eucal}
\usepackage{hyperref}

\newtheorem{theorem}{Theorem}[section]
\newtheorem{proposition}[theorem]{Proposition}
\newtheorem{corollary}[theorem]{Corollary}
\newtheorem{lemma}[theorem]{Lemma}
\newtheorem*{theorem*}{Theorem}
\newtheorem*{proposition*}{Proposition}
\newtheorem*{corollary*}{Corollary}
\newtheorem*{lemma*}{Lemma}
\theoremstyle{definition}
\newtheorem{definition}[theorem]{Definition}

\newtheorem{example}[theorem]{Example}
\newtheorem{remark}[theorem]{Remark}
\newtheorem*{remark*}{Remark}

\newtheorem*{definition*}{Definition}

\newcommand{\coring}[1]{\mathfrak{#1}}
\newcommand{\tensor}[1]{\otimes_{#1}}

\newcommand{\rcomod}[1]{ \mathsf{Comod}_{#1}}
\newcommand{\rmod}[1]{\mathsf{Mod}_{#1}}

\renewcommand{\hom}[3]{\mathrm{Hom}_{#1}(#2,\,#3)}

\newcommand{\rend}[2]{\mathrm{End}({#1}_{#2})}
\newcommand{\lend}[2]{\mathrm{End}({}_{#1}#2)}

\newcommand{\can}[1]{\mathsf{can}_{#1}}

\newcommand{\aut}[1]{\mathbf{Aut}(#1)}

\newcommand{\gl}[1]{\mathbf{Gr}(#1)}
\newcommand{\gla}[1]{\overline{\mathbf{Gr}}(#1)}
\newcommand{\galois}[1]{\mathbf{Gal}(#1)}
\newcommand{\swe}[2]{#1 \tensor{#2} #1}

\newcommand{\fk}[1]{\mathfrak{#1}}

\newcommand{\lr}[1]{\left(\underset{}{} #1 \right)}

\newcommand{\Scr}[1]{\mathscr{#1}}

\newcommand{\Sf}[1]{\mathsf{#1}}
\newcommand{\D}{\Scr{D}}
\newcommand{\HH}{\Scr{H}}

\begin{document}
\title{On the set of grouplikes of a coring\footnote{Research supported by  grant MTM2007-61673 from the Ministerio de
Educaci\'{o}n y Ciencia of Spain, and P06-FQM-01889 from Junta de Andaluc{\'\i}a}}
\author{L. El Kaoutit \\
\normalsize Departamento de \'{A}lgebra \\ \normalsize Facultad de Educaci\'{o}n y Humanidades \\
\normalsize  Universidad de Granada \\ \normalsize El Greco N$^o$
10, E-51002 Ceuta, Espa\~{n}a \\ \normalsize
e-mail:\textsf{kaoutit@ugr.es} \and J. G\'omez-Torrecillas \\
\normalsize Departamento de \'{A}lgebra \\ \normalsize Facultad de Ciencias \\
\normalsize Universidad
de Granada\\ \normalsize E18071 Granada, Espa\~{n}a \\
\normalsize e-mail: \textsf{gomezj@ugr.es} }

\date{\today}

\maketitle

\begin{abstract}
We focus our attention to the set $\gl{\coring{C}}$ of grouplike
elements of a coring $\coring{C}$ over a ring $A$. We do some
observations on the actions of the groups $U(A)$ and
$\aut{\coring{C}}$ of units of $A$ and of automorphisms of corings
of $\coring{C}$, respectively, on $\gl{\coring{C}}$, and on the
subset $\galois{\coring{C}}$ of all Galois grouplike elements.
Among them, we give conditions on $\coring{C}$ under which
$\galois{\coring{C}}$ is a group, in such a way that there is an
exact sequence of groups $\{1\} \rightarrow U(A^{g}) \rightarrow U(A) \rightarrow
\galois{\coring{C}} \rightarrow \{1\},$
where $A^g$ is the subalgebra of coinvariants for some $g \in
\galois{\coring{C}}$.
\end{abstract}

\section*{Introduction}
The concept of non-abelian cohomology of groups has been extended
to the framework of Hopf algebras by P. Nuss and M. Wambst in
\cite{Nuss/Wambst:2007, Nuss/Wambst:2008}. Given a Hopf algebra $H$, a right
$H$--comodule algebra $A$, and a right Hopf $H-A$--module $M$, the
first descent cohomology set $\D^1(H,M)$ of $H$ with coefficients
in $M$ is defined in terms of all Hopf module structures on $M$.
When $B \subseteq A$ is a $G$--Galois extension, where $G$ is a
finite group acting on $A$ by automorphisms, then by
\cite[Proposition 2.5]{Nuss/Wambst:2007} there is an isomorphism
of pointed sets $\D^1(K^G,M) \cong \HH^1(G,\aut{M_A})$, where the
last stands for the first non-abelian cohomology set of $G$ with
coefficients in $\aut{M_A}$ \cite{Serre:1950}. Here, $K^G$ is the Hopf algebra of
functions on the group $G$ in a commutative base ring $K$. In
\cite{Brzezinski:2008}, T. Brzezi\'nski has shown that this descent
cohomology can be satisfactorily extended to the framework of
comodules over corings, introducing the first descent cohomology
set $\D^1(\coring{C},M)$, where $\coring{C}$ is a coring over a
ring $A$, and $M$ is a right $\coring{C}$--comodule. Since the
definition of descent cohomology of \cite{Nuss/Wambst:2007} is a
special case of \cite[Definition 2.2]{Brzezinski:2008}, we know
that there must be an interpretation of the aforementioned
non-abelian cohomology set $\HH^1(G,\aut{M_A})$ in terms of descent
cohomology of a coring with coefficients in a comodule. The first
remark in this note (Theorem \ref{clasico}) gives such an
interpretation, in the case $M = A$. Our approach uses the fact
that the comodule structures on $A$ over an $A$--coring
$\coring{C}$ are parametrized by the set $\gl{\coring{C}}$ of all
grouplike elements of $\coring{C}$ \cite{Brzezinski:2002}. We thus
consider the particular case of \cite[Definition
2.2]{Brzezinski:2008} of the first descent cohomology set
$\D^1(\coring{C},g)$ of the $A$--coring $\coring{C}$ at a grouplike
element $g \in \gl{\coring{C}}$ (Definition \ref{atg}).

We focus our attention to the set $\gl{\coring{C}}$. We do some
observations on the actions of the groups $U(A)$ and
$\aut{\coring{C}}$ of units of $A$ and of automorphisms of corings
of $\coring{C}$, respectively, on $\gl{\coring{C}}$. Among them,
let us mention that if $\D^{1}(\coring{C},g) = \{1 \}$, then
$\aut{\coring{C}}$ is isomorphic to a quotient group
$U(A)_g/U(A^g)$, where $A^g$ (Definition \ref{Def-1}) is the subring of $g$--coinvariants of
$A$, and $U(A)_g$ is a subgroup of $U(A)$ (Corollary \ref{obs1}).
We also give (Theorem \ref{mejor}) conditions under which the set
$\galois{\coring{C}}$ of all Galois grouplike elements is a group,
in such a way that there is an exact sequence of groups
\[
\{1\} \rightarrow U(A^{g}) \rightarrow U(A) \rightarrow
\galois{\coring{C}} \rightarrow \{1\}.
\]
We also give some conditions on $g \in \galois{\coring{C}}$ to
have that $\D^1(\coring{C},g) = \{1 \}$. Our approach here makes
use of the theory of cosemisimple corings developed in
\cite{ElKaoutit/Gomez/Lobillo:2004c} and
\cite{ElKaoutit/Gomez:2003a}.

Some examples illustrate our results.

\section{Grouplikes, non abelian cohomology and descent cohomology}\label{Sec-1}

Let $(\coring{C},\Delta_{\coring{C}}, \varepsilon_{\coring{C}})$
be a coring over a $K$--algebra $A$ ($K$ is a commutative ring).
Thus, $\coring{C}$ is an $A$--bimodule, and $\Delta_\coring{C} :
\coring{C} \rightarrow \coring{C} \tensor{A} \coring{C}$ and
$\varepsilon_\coring{C} : \coring{C} \rightarrow A$ are
homomorphisms of $A$--bimodules subject to axioms of
coassociativity and counitality:
$(\coring{C}\tensor{A}\Delta_{\coring{C}}) \circ
\Delta_{\coring{C}} = (\Delta_{\coring{C}}\tensor{A}\coring{C})
\circ \Delta_{\coring{C}}$ and
$(\coring{C}\tensor{A}\varepsilon_{\coring{C}}) \circ
\Delta_{\coring{C}} = (\varepsilon_{\coring{C}}
\tensor{A}\coring{C}) \circ \Delta_{\coring{C}} = \coring{C}$. A
right $\coring{C}$--comodule is a pair $(M,\rho_{M})$ consisting
of a right $A$--module and a right $A$--linear map $\rho_{M}: M
\rightarrow M\tensor{A}\coring{C}$, called right
$\coring{C}$--coaction, such that
$(M\tensor{A}\Delta_{\coring{C}}) \circ \rho_M =
(\rho_M\tensor{A}\coring{C}) \circ \rho_M$ and
$(M\tensor{A}\varepsilon_{\coring{C}}) \circ \rho_M=M$. A morphism
of right $\coring{C}$--comodules $f : (M,\rho_M) \rightarrow
(N,\rho_N)$ is a right $A$--linear map $f :M \rightarrow N$ such
that $ \rho_N  \circ  f = (f \tensor{A} \coring{C})  \circ
\rho_M$. With these morphisms, right $\coring{C}$--comodules form
a category. Details on corings and their comodules are easily
available in \cite{Brzezinski/Wisbauer:2003}.

\begin{definition}\label{Def-1}
An element $g \in\coring{C}$ is said to be a \emph{grouplike
element} if $\Delta_{\coring{C}}(g) = g \tensor{A}g $ and
$\varepsilon_{\coring{C}}(g)=1$. The set of all grouplike elements
of $\coring{C}$ will be denoted by $\gl{\coring{C}}$. The subring of $g$--\emph{coinvariant} elements is defined by $$A^g=\{a \in
A|\quad ag \, =\, ga \}.$$
\end{definition}

\begin{example}\label{Sweedler-corg}
If $B \rightarrow A$ is any ring extension, and $A\tensor{B}A$ is
its associated Sweedler's $A$--coring with comultiplication
$a\tensor{B}a' \mapsto (a\tensor{B}1)\tensor{A}(1\tensor{B}a')$
and counit the multiplication map $a\tensor{B}a' \mapsto aa'$,
$a,a'\in A$, then it is clear that $1\tensor{B}1 \in
\gl{A\tensor{B}A}$. Given an SBN (Single Basis Number) ring $A$,
then by \cite[p. 113]{Anderson/Fuller:1974}, there exist elements
$a,a',b,b'\in A$ such that $$ ab + a'b'\,=\, 1,\quad
ba\,=\,b'a'\,=\,1,\quad\text{ and }\,\, b'a\,=\,ba'\,=\,0. $$
Clearly $g\,=\, a\tensor{\mathbb{Z}}b + a'\tensor{\mathbb{Z}}b'$
is a grouplike element of the Sweedler $A$--coring
$A\tensor{\mathbb{Z}}A$.
\end{example}

\begin{example}
Consider a coring $\coring{C}$ (resp. $\coring{C}'$) over a ring
$A$ (resp. $A'$). Let $\rho : A \rightarrow A'$ be a homomorphism
of rings, and consider a homomorphism of corings $\varphi :
\coring{C} \rightarrow \coring{C}'$ in the sense of
\cite{Gomez:2002}. This morphism restricts to a map $\varphi:
\gl{\coring{C}} \rightarrow \gl{\coring{D}}$. Moreover, for each
$g \in \gl{\coring{C}}$, $\rho$ induces a homomorphism of rings
$\rho: A^{g} \rightarrow A'^{\varphi(g)}$.
\end{example}

It is known \cite[Lemma 5.1]{Brzezinski:2002} that there is a
bijection between $\gl{\coring{C}}$ and the set of all right
$\coring{C}$-coactions on the right  module $A_A$. Let $[g]A$
denote the right $\coring{C}$--comodule structure defined on $A$
by $g \in \gl{\coring{C}}$. The right $\coring{C}$-coaction
$\rho_{[g]A}: [g]A \rightarrow [g]A\tensor{A}\coring{C} \cong
\coring{C}$ is given by sending $a \mapsto ga$.
 Conversely any right $\coring{C}$-coaction $\rho_A$ determines a
unique element $\rho_A(1) \in \gl{\coring{C}}$. A similar
bijection exists taking left $\coring{C}$-coactions on the left
module ${}_AA$. We denote by $A[g]$ the left $\coring{C}$-comodule
induced by $g \in \gl{\coring{C}}$.

It is easily checked that the subring $A^g$ of $A$ can be identified with both rings of
endomorphisms of the right $\coring{C}$--comodule $[g]A$ and of
the left $\coring{C}$--comodule $A[g]$. That is, $A^g =
\rend{[g]A}{\coring{C}} = \lend{\coring{C}}{A[g]}$. In fact, for
two grouplike elements $g,h \in \gl{\coring{C}}$, we have
$$\hom{\coring{C}}{[g]A}{[h]A} \, = \, \left\{\underset{}{} \alpha \in A|\,\, \alpha
g = h \alpha \right\}.$$ Therefore, $[g]A\,\cong\,[h]A$ as right
$\coring{C}$-comodules if and only if $g$ and $h$ are
\emph{conjugated} in the sense that there exists $\alpha \in U(A)$
such that $h\,=\, \alpha g \alpha^{-1}$. These remarks suggest, in
view of \cite{Brzezinski:2008}, the following definition, due to
T. Brzezi\'nski.

\begin{definition}\label{atg}
Consider the action of the group of units $U(A)$ of $A$ on
$\gl{\coring{C}}$
\begin{equation}\label{accion1}
\xymatrix@R=0pt@C=50pt{ U(A)
\,\times \, \gl{\coring{C}} \ar@{->}[r] & \gl{\coring{C}} \\
(\alpha, g) \ar@{|->}[r] & \alpha g \alpha^{-1}.}
\end{equation}
Let $\gla{\coring{C}}$ denote the quotient set of
$\gl{\coring{C}}$ under the action \eqref{accion1}. If
$\gl{\coring{C}}$ is not empty, then for each $g \in
\gl{\coring{C}}$ we can define the \emph{pointed set of descent
$1$--cocycles on $\coring{C}$ with coefficients in $[g]A$} as
\[
Z^{1}(\coring{C},[g]A) := (\gl{\coring{C}},g),
\]
and the \emph{first cohomology pointed set of $\coring{C}$ with
coefficients in $[g]A$} as
\[
\D^{1}(\coring{C},[g]A) := (\gla{\coring{C}},\overline{g}),
\]
where $(X,x)$ means a pointed set with a distinguished element $ x \in X$. We shall use the simplified notations $Z^{1}(\coring{C},g)$ and
$\D^{1}(\coring{C},g)$, respectively, and we will refer to them as the
\emph{pointed set of descent $1$--cocycles on $\coring{C}$ at
$g$}, and the \emph{first descent cohomology of $\coring{C}$ at
$g$}, respectively.
The \emph{zeroth descent cohomology group} of $\coring{C}$ at $g$ is defined to
be the group of $\coring{C}$-comodules automorphisms of $[g]A$, and can be identified with the group of units $U(A^g)$ of the ring $A^g$, i.e., $$\D^0(\coring{C},g)\,=\, U(A^g).$$
\end{definition}

Our first aim is to exhibit a direct evidence of the fact that
$\D^{1}(\coring{C},g)$ is a genuine version for corings of Serre's
nonabelian cohomology of groups.

\begin{example}\label{cocyclos}
Let $G$ be a finite group acting by automorphisms on a ring $A$.
Consider $R\,=\, G*A$ the associated crossed product. As $R$ is a
free right $A$--module with basis $G$, its right dual
$R^*\,=\,\hom{A}{R}{A}$ is an  $A$-coring according to
\cite[Theorem 3.7]{Sweedler:1975} (with comultiplication and
counit induced by the duals of the multiplication and the unit of
the $A$-ring $R$). Our next aim is to establish a bijection
between $\gl{R^*}$ and the set of all non abelian $1$-cocycles
$Z^1(G^{op},\,U(A))$ in the sense of \cite{Serre:1950}. Of
course here the action of the opposite group $G^{op}$ on the group
$U(A)$ is induced by the given action of the group $G$ on the ring
$A$. So we denote this action by $\alpha^x$ for every $\alpha \in
U(A)$ and $x \in G^{op}$.

\begin{proposition}\label{cociclosSerre}
The map $\Theta : \gl{R^*} \rightarrow Z^1(G^{op},U(A))$
which sends $h \in \gl{R^*}$ to its restriction to $G$ is a
bijection. Under this bijection, the trivial $1$--cocycle
corresponds to the grouplike given by the trace map $\fk{t} : R
\rightarrow A$ defined by $\fk{t}(\sum_{x\in G} x a_x) = \sum_{x
\in G}a_x$. Moreover, $A^{\fk{t}}$ coincides with the subring of
the $G$--invariant elements of $A$.
\end{proposition}
\begin{proof}
Let us denote by $\{x,x^*\}_{x \in\, G} \subseteq R \times R^*$
the finite dual basis of the right free $A$-module $R_A$ given by
$G$.  The comultiplication and the counit of the $A$-coring $R^*$
are defined as follows
$$\xymatrix@R=0pt@C=40pt{ R^* \ar@{->}^-{\Delta}[r] &  R^*\tensor{A}R^* \\
\varphi \ar@{|->}[r] & \sum_{x\in \, G}\varphi x \tensor{A}x^*,  } \qquad
\xymatrix@R=0pt@C=40pt{ R^* \ar@{->}^-{\varepsilon}[r] &  A \\ \varphi \ar@{|->}[r] &
\varphi(1_R). }$$
We have an isomorphism
$$\xymatrix{\Upsilon: R^*\tensor{A}R^* \ar@{->}[r] &
(R\tensor{A}R)^*}, \quad \lr{\varphi\tensor{A}\psi \longmapsto
\left[\underset{}{} r\tensor{A}t \mapsto \varphi(\psi(r)t)
\right]} $$ Now, a right $A$-linear map $ h : R \to A$ belongs to
$\gl{R^*}$ if and only if
\begin{equation}\label{Eq-grouplike}
h(1_R)\,=\,1_A \,\,\text{ and }\,\,
\sum_{x\in\,G}h\,x\tensor{A}x^*\,\,=\,\, h\tensor{A}h.
\end{equation}
So given $h \in \gl{R^*}$, and applying $\Upsilon$ to the second
equality in \eqref{Eq-grouplike}, we obtain the  equality  $h(x
y)\,=\, h(y)\,h(x)^y$, for every pair of elements $x,y \in G$.
Taking $y=x^{-1}$, we get $h(x^{-1}) h(x)^{x^{-1}}= 1_A  =  h(x)
h(x^{-1})^{x}$, since $h(1_R)=1_A$. Applying $x$ to the
equality $h(x^{-1}) h(x)^{x^{-1}}= 1_A$, we obtain $h(x)h(x^{-1})^{x}=
h(x^{-1})^{x}h(x)=1_A$. That is, $h(x) \in U(A)$, for every $x \in
G$. Concluding, we have defined a map
\begin{equation}\label{tilda}
\xymatrix@R=0pt{ \gl{R^*} \ar@{->}^-{\Theta}[rr] & &
Z^1(G^{op},U(A)) \\ h \ar@{|->}[rr] & &
\left[\underset{}{}\Theta(h): x \longmapsto h(x)  \right].}
\end{equation}
It is clear that $\Theta$ is injective, since $G$ is a basis for
the right $A$-module $R$. Let us check that it is also surjective.
Consider  any $1$-cocycle $f:G^{op} \to U(A)$, and define
$\widehat{f}: R \to A$ by sending  $x * a \mapsto f(x)a$ for $x
\in G$ and $a \in A$. Clearly $\widehat{f}$ is a right $A$-linear
map, and $\varepsilon(\widehat{f})=\widehat{f}(1_R)= f(\Sf{e})1_A
= f(\Sf{e})$ (here $\Sf{e}$ is the neutral element of $G$). By the
$1$-cocycle condition on $f$, we know that
$f(\Sf{e})=f(\Sf{e})^2$, that is, $f(\Sf{e})=1_A$ and so
$\varepsilon(\widehat{f})=1_A$. Now an easy computation using
again the $1$-cocycle condition shows that
$$\Upsilon\lr{\sum_{x\in\,G}
\widehat{f}x\tensor{A}x^*}(y\tensor{A}z)\,\,=\,\,
\Upsilon\lr{\widehat{f}\tensor{A}\widehat{f}}(y\tensor{A}z),$$ for
every pair of elements $y,z \in G$, which implies that
$\Delta(\widehat{f})\,=\,\widehat{f}\tensor{A}\widehat{f}$.
Therefore, $\widehat{f} \in \gl{R^*}$. Obviously, we have
$\Theta(\widehat{f}) =f$, and this establishes the desired
surjectivety. Clearly the distinguished $1$-cocycle $\fk{e}:
G^{op} \to U(A)$ sending $x \mapsto 1$ corresponds  then to the
grouplike element $\fk{t} : R \to A$ defined by $\sum_{x\in G} x
a_x \mapsto \sum_{x \in G}a_x$.  The coinvariant ring $A^{\fk{t}}$
coincides with the invariant subring of $A$ with respect to the
$G$-action i.e. $A^{\fk{t}}\,=\,\{a \in A|\, x(a)\,=\,a, \forall x
\in G\}$.
\end{proof}
\end{example}

Recall from \cite{Serre:1950}, that two $1$-cocycles $f$ and $h$
are cohomologous if there exists $\alpha \in U(A)$ such that
$f(x)\,=\, \alpha^{-1} h(x) \alpha^x$, for every $x \in G^{op}$.
Using the bijection \eqref{tilda}, we can easily check that two
$1$-cocycles are cohomologous if and only if their corresponding
grouplike elements are conjugated. On the other hand, the equality $A^{\fk{t}}=A^{G}$ clearly implies that $U(A^{\fk{t}}) = U(A)^{G^{op}} $, where the later stands for $\HH^0(G^{op},U(A))$ the zeroth non-abelian cohomology group as in \cite{Serre:1950}.
Therefore, we deduce from Proposition \ref{cociclosSerre}:

\begin{theorem}\label{clasico}
The map $\Theta$ of Proposition \ref{cociclosSerre} induces an
isomorphism of pointed sets
\[
\D^1(R^*,\mathfrak{t}) \cong \HH^1(G^{op},U(A)),
\]
and there is an equality of groups $$\D^0(R^*, \fk{t})\,=\, \HH^0(G^{op},U(A)).$$
\end{theorem}

\begin{remark}
Since the coring $R^*$ is finitely generated and projective as a
left $A$-module, its category of right comodules is isomorphic to
the category of right $R$-modules. Taking this into account, one
can adapt the proof of Proposition \ref{cociclosSerre} in order to
show that for every right $A^{\mathfrak{t}}$-module $N$, there is
an isomorphism of pointed sets
$$ \D^1(R^*, N\tensor{A^{\mathfrak{t}}}[\mathfrak{t}]A) \,\cong \, \HH^1(G^{op},{\bf Aut}_A(N\tensor{A^{\mathfrak{t}}}A)),$$
and an equality of groups $$ \D^0(R^*,
N\tensor{A^{\mathfrak{t}}}[\mathfrak{t}]A) \,=\, \HH^0(G^{op},{\bf
Aut}_A(N\tensor{A^{\Sf{t}}}A)),$$ where for every right
$\coring{C}$-comodule $M$, $\D^{\bullet}(\coring{C},M)$ are
defined as in \cite{Brzezinski:2008}, and ${\bf Aut}_A(M)$ is the
group of all automorphisms of the underlying right $A$--module of
$M$.
\end{remark}

\section{Groups acting on grouplikes}

The maps defined in the following lemma will be used in the sequel
where the role of the extension $B \rightarrow A$ will be played
by the inclusions $A^g \subseteq A$, and where $g$ runs
$\gl{\coring{C}}$ whenever $\gl{\coring{C}} \neq \emptyset$.

\begin{lemma}\label{auto}
Let $B \rightarrow A$ be any ring extension.
\begin{enumerate}[(a)]
\item Let $\alpha \in U(A)$ and consider the subring
    $\alpha^{-1} B \alpha$ of $A$. Then the map
$$\xymatrix@R=0pt@C=50pt{\psi_{\alpha}: A\tensor{B} A
\ar@{->}[r] & A\tensor{\alpha^{-1} B \alpha}A
\\ a\tensor{B}a' \ar@{|->}[r] & a\alpha\tensor{\alpha^{-1} B\alpha}\alpha^{-1} a'
}$$ is an isomorphism of $A$--corings.

\item The map $$\psi_{-}: \left\{\underset{}{}\alpha \in
    U(A)|\, \alpha^{-1}B\alpha \,=\, B\right\} \longrightarrow
    \aut{A\tensor{B}A}$$ defines an anti-homomorphism of
    groups.
\end{enumerate}
\end{lemma}
\begin{proof}$(a)$ We only prove that $\psi_{\alpha}$
is a well defined map. So, for every $a,a' \in A$ and $b\in B$, we
have
\begin{eqnarray*}
  \psi_{\alpha}(ab\tensor{K}a') &=& ab\alpha \tensor{\alpha^{-1} B \alpha} \alpha^{-1} a' \\
   &=& a\alpha(\alpha^{-1}b\alpha)\tensor{\alpha^{-1} B \alpha} a' \\
   &=& a\alpha \tensor{\alpha^{-1} B\alpha} (\alpha^{-1}b\alpha)\alpha^{-1}a' \\ &=&
   a\alpha\tensor{\alpha^{-1}B\alpha}\alpha^{-1}ba' \, = \,
   \psi_{\alpha}(a\tensor{K}ba'),
\end{eqnarray*} that is, $\psi_{\alpha}$ is a well defined map.

$(b)$ Straightforward.
\end{proof}

Every grouplike element $g \in \gl{\coring{C}}$ of the $A$--coring
$\coring{C}$ defines a \emph{canonical morphism} of $A$--corings:
$$\can{g}: \swe{A}{A^g} \longrightarrow \coring{C},\quad
\lr{a\tensor{A^g}a' \longmapsto aga'}.$$

On the other hand, a straightforward computation shows that
\begin{equation}\label{iguales}
A^{\alpha g \alpha^{-1}} = \alpha A^g \alpha^{-1} \quad \hbox{ for all } \alpha \in U(A).
\end{equation}

Moreover, for every $\alpha \in U(A)$, we have the commutative
diagram of homomorphisms of $A$--corings:
\begin{equation}\label{triangulo}
\xymatrix@C=50pt{ \swe{A}{A^{\alpha g \alpha^{-1}}}
\ar@{->}^-{\can{\alpha g
\alpha^{-1}}}[rd]  & \\ & \coring{C} \\
\swe{A}{A^g}. \ar@{->}_-{\can{g}}[ru] \ar@{->}_-{\cong}^-{\psi_{\alpha^{-1}}}[uu] & }
\end{equation}

Recall from \cite{Brzezinski:2002} that a grouplike $g \in
\gl{\coring{C}}$ is said to be \emph{Galois} if $\can{g}$ is
bijective. It follows from diagram \eqref{triangulo} that $g$ is
Galois if and only if $\alpha g \alpha^{-1}$ is Galois. Thus, if
we denote by $\galois{\coring{C}}$ the set of all Galois grouplike
elements of $\coring{C}$, then the action \eqref{accion1} restricts to an action
\begin{equation*}
\xymatrix@R=0pt@C=50pt{ U(A)
\,\times \, \galois{\coring{C}} \ar@{->}[r] & \galois{\coring{C}} \\
(\alpha, g) \ar@{|->}[r] & \alpha g \alpha^{-1}.}
\end{equation*}

The group $\aut{\coring{C}}$ of all $A$--coring automorphisms of
$\coring{C}$ acts obviously on $\gl{\coring{C}}$:
\begin{equation}\label{accion2}
\xymatrix@R=0pt@C=50pt{ \aut{\coring{C}} \,
\times \, \gl{\coring{C}} \ar@{->}[r] & \gl{\coring{C}} \\
(\varphi, g) \ar@{|->}[r] & \varphi \cdot g := \varphi(g). }
\end{equation}
Since every $\varphi \in \aut{\coring{C}}$ is, in particular, a
homomorphism of $A$--bimodules, it follows that the actions
\eqref{accion2} and \eqref{accion1} commute, that is
\begin{equation*}
\varphi \cdot (\alpha \cdot g) = \alpha \cdot (\varphi \cdot g), \quad \forall g \in \gl{\coring{C}},
\forall \alpha \in U(A), \forall \varphi \in \aut{\coring{C}}.
\end{equation*}

The action \eqref{accion2} restricts to an action
\begin{equation}\label{Gaccion2}
\xymatrix@R=0pt@C=50pt{ \aut{\coring{C}} \,
\times \, \galois{\coring{C}} \ar@{->}[r] & \galois{\coring{C}} \\
(\varphi, g) \ar@{|->}[r] &  \varphi(g), }
\end{equation}
as the following proposition shows.

\begin{proposition}\label{auto-accion}
\begin{enumerate}[(1)]
\item\label{autcor} For every element $g \in \gl{\coring{C}}$
    and $\varphi \in \aut{\coring{C}}$, we have $ A^g \,=\,
    A^{\varphi(g)}$. Moreover, the following diagram of
    morphisms of $A$--coring commutes
$$\xymatrix@C=50pt{ \swe{A}{A^g}
\ar@{->}^-{\can{g}}[r] \ar@{=}[d] & \coring{C} \ar@{->}^-{\varphi}[d] \\
\swe{A}{A^{\varphi(g)}} \ar@{->}^-{\can{\varphi(g)}}[r] &
\coring{C}. }$$ Therefore,  $g \in \galois{\coring{C}}$ if and
only if $\varphi(g) \in \galois{\coring{C}}$.
\item\label{autcorgal} If $g, h \in \galois{\coring{C}}$, then
    $A^g = A^h$ if and only if there exists $\varphi \in
    \aut{\coring{C}}$ such that $\varphi(h) = g$.
\end{enumerate}
\end{proposition}
\begin{proof}
\eqref{autcor} We only prove that $A^g = A^{\varphi(g)}$. Start
with an element $b \in A^g$, then $ b \varphi(g) = \varphi(bg) =
\varphi(gb) =\varphi(g) b$ because $\varphi$ is a homomorphism of
$A$--bimodules, which implies that $b \in A^{\varphi(g)}$. Thus
$$ A^g \, \subseteq \, A^{\varphi(g)} \, \subseteq
A^{\varphi^{-1}(\varphi(g))}\, = \, A^g.$$
\eqref{autcorgal} Assume that $A^g = A^h$ for $g,h \in
\galois{\coring{C}}$. Then $\varphi = \can{g} \circ \can{h}^{-1}
\in \aut{\coring{C}}$ and $\varphi (h) = g$. The converse follows from item $(1)$.
\end{proof}

For every element $g \in \gl{\coring{C}}$, we define
$$U(A)_g\,=\, \left\{\underset{}{} \alpha \in U(A)|\,\, \alpha A^g =
A^g \alpha
\right\}\, = \, \left\{\underset{}{} \alpha \in U(A)|\,\,
A^{\alpha g \alpha^{-1}} = A^g \right\},$$ where in the second
equality, we have used equation \eqref{iguales}. It is clear that
$U(A)_g$ is a subgroup of $U(A)$ which contains the group of units
$U(A^g)$ of the subring $A^g$.

\begin{proposition}\label{conju-sucesion}
Let $\coring{C}$ be an $A$--coring.
\begin{enumerate}[(a)]
\item\label{suca} For every element $g \in \gl{\coring{C}}$,
    $U(A^g)$ is a normal subgroup of $U(A)_g$.
\item \label{sucb} For every $g \in \gl{\coring{C}}$ and
    $\beta \in U(A)$, we have $$ \beta U(A)_g \beta^{-1} \, =
    \, U(A)_{\beta g \beta^{-1}} .
$$
\item\label{succ} If $g \in \galois{\coring{C}}$, then there exists an
    exact sequence of groups $$\xymatrix@R=0pt{1 \ar@{->}[r] &
    U(A^g) \ar@{->}[r] & U(A)_g \ar@{->}^-{\phi_g}[r] &
    \aut{\coring{C}} \\ & & \alpha \ar@{|->}[r] & \can{g}
    \circ \psi_{\alpha^{-1}} \circ \can{g}^{-1}. }$$
\item\label{sucd} If $\galois{\coring{C}}$ is non empty and
    the action of $U(A)$ on $\galois{\coring{C}}$ is
    transitive, then, for every $g \in \galois{\coring{C}}$,
    $\phi_g$ is surjective and, thus, we have an isomorphism
    of groups
    \begin{equation*}
\aut{\coring{C}} \cong U(A)_g/U(A^g).
    \end{equation*}
\end{enumerate}
\end{proposition}
\begin{proof}
\eqref{suca} Let $\beta $ be an arbitrary element in $U(A)_g$.
Given an element $\alpha \in U(A^g)$, by definition there exists
$\gamma \in A^g$ such that $\beta \alpha \beta^{-1}\,=\, \gamma$, and so
$\gamma \in U(A^g)$. Therefore, $\beta U(A^g) \beta^{-1}
\subseteq U(A^g)$.

\eqref{sucb} It follows from the fact that $U(A)_g$ is the
stabilizer in $U(A)$ of $A^g$ for the action by conjugation of
$U(A)$ on the set of all subalgebras of $A$.

\eqref{succ} An element $\alpha \in U(A)_g$ is such that
$\phi_g(\alpha) =1$ if and only if $\can{g} \circ
\psi_{\alpha^{-1}} = \can{g}$ if and only if $\can{g} \circ
\psi_{\alpha^{-1}} (1 \tensor{A^g} 1) = \can{g} (1 \tensor{A^g}
1)$ if and only if $\alpha^{-1} g \alpha = g$ if and only if
$\alpha \in U(A^g)$, and the exactness follows.

\eqref{sucd} Let $g \in \galois{\coring{C}}$ and $\varphi \in
\aut{\coring{C}}$. Obviously, $\varphi(g) \in \galois{\coring{C}}$
and, since $\galois{\coring{C}} = \{ \beta g \beta^{-1} : \beta
\in U(A) \}$, there exists $\alpha \in U(A)$ such that $\varphi(g)
= \alpha^{-1} g \alpha$. We know that $$A^g = A^{\varphi(g)} =
A^{\alpha^{-1} g \alpha} = \alpha^{-1}A^g\alpha,$$ that is, $\alpha
\in U(A)_g$. Moreover, it is easily checked that
$\phi_g(\alpha)(g) = \alpha^{-1}g\alpha$ and, since $g$ generates
$\coring{C}$ as an $A$--bimodule, this implies that $\varphi =
\phi_g(\alpha)$. Therefore, $\phi_g$ is surjective.
\end{proof}

\begin{corollary}\label{obs1}
If $g$ is a Galois grouplike of $\coring{C}$ such that
$\D^1(\coring{C},g) = \{1\}$, then $$\aut{\coring{C}} \cong
U(A)_g/U(A^g).$$
\end{corollary}

\begin{remark}
When $A$ is commutative, Corollary \ref{obs1} says that
$\aut{\coring{C}}$ is the \emph{coGalois group} of the extension
$A^g \subseteq A$, see \cite{Masuoka:1989b} for the case of field
extensions.
\end{remark}

\begin{theorem}\label{mejor}
Let $\coring{C}$ be an $A$--coring such that there exists $g \in
\galois{\coring{C}}$ and the action of $U(A)$ on
$\galois{\coring{C}}$ is transitive (e.g. $\D^1(\coring{C},g) = \{
1 \}$). The following statements are equivalent:
\begin{enumerate}[(i)]
\item\label{trans1} $U(A)_g = U(A)$ (i.e., $\alpha A^g = A^g
    \alpha$ for every $\alpha \in U(A)$);
\item\label{trans2} $U(A)_h = U(A)$ for every $h \in
    \galois{\coring{C}}$;
\item\label{trans3} the action of $\aut{\coring{C}}$ on
    $\galois{\coring{C}}$ is transitive.
\end{enumerate}
Furthermore, if one of these equivalent conditions is satisfied,
then $A^h = A^g$ for every $h \in \galois{\coring{C}}$, and the
map $\xi_g : \aut{\coring{C}} \rightarrow
        \galois{\coring{C}} $ defined by $\xi_g (\varphi) =
        \varphi (g)$ for $\varphi \in \aut{\coring{C}}$ is
        bijective and, thus, $\galois{\coring{C}}$ can be endowed
        with the structure of a group. Moreover, there exists
        a short exact sequence of groups
\[
\{1\} \rightarrow U(A^{g}) \rightarrow U(A) \rightarrow
\galois{\coring{C}} \rightarrow \{1\}.
\]
\end{theorem}
\begin{proof}
\eqref{trans1} $\Rightarrow$ \eqref{trans3}. By assumption we have
$U(A)_g = U(A)$. On the other hand, every grouplike is of the form
$\alpha g \alpha^{-1}$ for some $\alpha \in U(A)$. Now, $\alpha g
\alpha^{-1} = \phi_g(\alpha^{-1})(g)$, where $\phi_g(\alpha^{-1}) \in
\aut{\coring{C}}$ is given by Proposition
\ref{conju-sucesion}.\eqref{succ}. This means that each grouplike element is in the obrit of $g$ under the action \eqref{accion2}. \\
\eqref{trans3} $\Rightarrow$ \eqref{trans2}. Given $h \in
\galois{\coring{C}}$ and $\alpha \in U(A)$, we know from
\eqref{triangulo} that $\alpha h \alpha^{-1} \in
\galois{\coring{C}}$. Since the action of $\aut{\coring{C}}$ on
$\galois{\coring{C}}$ is transitive, there is $\varphi \in
\aut{\coring{C}}$ such that $\varphi (h) = \alpha h \alpha^{-1}$.
Proposition \ref{conju-sucesion}.(2) and equation \eqref{iguales}
give now that
\[
A^h = A^{\varphi(h)} = A^{\alpha h \alpha^{-1}} = \alpha
A^h \alpha^{-1},
\]
that is, $\alpha \in U(A)_h$.
\\ Since \eqref{trans2} $\Rightarrow$ \eqref{trans1} is obvious, the proof
of the equivalence between the three statements is done.

If $h \in \galois{\coring{C}}$, then $A^{h} = A^{\varphi(g)} = A^g$
for some $\varphi \in \aut{\coring{C}}$. On the other hand, it is
clear from assumption that $\xi_g$ is surjective. Since $g$, being Galois,
generates $\coring{C}$ as an $A$--bimodule, it follows that the
action of every automorphism of $\coring{C}$ on $g$ determines it
completely. Thus, $\xi_g$ is injective. Finally, the short exact
sequence of groups is given by Proposition
\ref{conju-sucesion}.\eqref{succ}-\eqref{sucd}.
\end{proof}

A ring $A$ is said to be a \emph{right invariant basis number
ring} (right \emph{IBN} ring for short), if $A^{(n)} \cong
A^{(m)}$ (direct sums of copies of $A$) as right $A$--modules for
 $n$, $m \in \mathbb{N}$ implies that $n=m$, see \cite[p.
114]{Anderson/Fuller:1974}. An $A$--coring $\coring{C}$ is said to
\emph{cosemisimple} if ${}_A\coring{C}$ is a flat module and every
right $\coring{C}$--comodule is semisimple, equivalently,
$\coring{C}_A$ is flat and every left $\coring{C}$--comodule is
semisimple. A \emph{simple cosemisimple} coring is a cosemisimple
coring with one type of simple right comodule or equivalently with
one type of simple left comodule; see \cite[Therorem
4.4]{ElKaoutit/Gomez:2003a} for a structure theorem of all
cosemisimple corings over an arbitrary ring. By \cite[Theorem
4.3]{ElKaoutit/Gomez/Lobillo:2004c}, every grouplike of a simple
cosemisimple coring $\coring{C}$ is Galois, that is,
$\gl{\coring{C}} = \galois{\coring{C}}$.

\begin{theorem}\label{simple-cosem}
Let $\coring{C}$ be an $A$--coring, and assume that there exists
$g \in \galois{\coring{C}}$. Assume that either $A^{g}$ is a
    division ring and $A$ is a right (or left) IBN ring, or $A$ is a division ring.  Then
    $$\gl{\coring{C}} = \galois{\coring{C}} = \left\{ \underset{}{} \alpha g \alpha^{-1}|\, \,
    \alpha \in U(A) \right\},$$ and, in particular,
    $\D^1(\coring{C},g) = \{ 1 \}$.
\end{theorem}
\begin{proof}
Assume first that $A^g$ is a division ring and $A$ is left or
right IBN. By \cite[Theorem 4.4]{ElKaoutit/Gomez/Lobillo:2004c}
(see also \cite[Theorem 3.10, Proposition
4.2]{ElKaoutit/Gomez:2003a}), $\coring{C}$ is a simple
cosemisimple $A$--coring and the functor $ - \tensor{A^g} [g]A :
\rmod{A^g} \rightarrow \rcomod{\coring{C}}$ is an equivalence of
categories. Thus, $[g]A \cong A^g \tensor{A^g} [g]A$ is a simple
right comodule. Given  $h \in \gl{\coring{C}}$, we have the right
$\coring{C}$--comodule $[h]A$. Since $\coring{C}$ is cosemisimple
with a unique type of simple right comodule represented by $[g]A$,
there is an isomorphism of right $\coring{C}$--comodules $[h]A
\cong \left( [g]A\right)^{(n)}$ (direct sum of copies of $[g]A$),
for some non zero natural number $n$. This isomorphism is, in
particular, an isomorphism of right free $A$--modules. Hence,
$n=1$ since $A$ is right IBN ring. Therefore, $[h]A \cong [g]A$,
as comodules, which means that $h = \alpha g \alpha^{-1}$ for some
$\alpha \in U(A)$, and we have done. In the case that $A$ is a
division ring, it is easy to show that $A^g$ is a division ring.
\end{proof}

\begin{corollary}\label{BA}
Let $B \subseteq A$ be a ring extension, and $A \tensor{B} A$ its
canonical Sweedler's coring.
\begin{enumerate}[(1)]
\item If $B$ is a division ring and $A$ is a right or left IBN
    ring, then $$\gl{\swe{A}{B}} = \left\{ \underset{}{}
    \alpha\tensor{B}\alpha^{-1}
|\,\, \alpha \in U(A) \right\}.$$
\item If $B \subseteq A$ is an extension of division rings and
    $\alpha B = B \alpha$ for every $\alpha \in A$, then
    $\gl{\swe{A}{B}}$ is a group isomorphic to $A^{\times}/B^{\times}$.
\end{enumerate}
\end{corollary}
\begin{proof}
By \cite[Proposition 4.2]{ElKaoutit/Gomez/Lobillo:2004c}, $B =
A^{1 \otimes_B 1}$. The corollary follows now from Theorem
\ref{simple-cosem}.
\end{proof}

\begin{example}
Let $G$ be a finite group acting on a division ring $A$ as in
Example \ref{cocyclos}, and let $T$ be the (division) subring of
all $G$--invariant elements of $A$. We know that $T =
A^{\mathfrak{t}}$. Assume that the trace map $\mathfrak{t}$ is a
Galois grouplike of $R^*$, where $R = G*A$. This means that $T
\subseteq A$ is Galois in the sense that the canonical map $G*A
\rightarrow \lend{T}{A}$ is bijective \cite{Kanzaki:1964}. Then, by Theorems
\ref{clasico} and \ref{simple-cosem}, $\HH^1(G^{op},A^{\times}) =
\{1\}$. This is a version of Hilbert's 90 theorem for division
rings.
\end{example}

\begin{remark}
The condition $\alpha B = B \alpha$ for every $\alpha \in A$ in
Corollary \ref{BA} is rather strong. An easy example is the
following. Let $A = \mathbb{C}_q(X,Y)$ the (noncommutative) field
of fractions of the complex quantum plane $\mathbb{C}_q[X,Y]$, and
$B = \mathbb{C}(X)$, the field of complex rational functions in
the variable $X$ (here, $q \neq 1$ is a complex number). It is
easy to show that $(1 + Y) B \neq B (1 + Y)$. In fact, $(1 + aY) B
\neq B (1 + aY)$ for infinitely many $a \in \mathbb{C}^{\times}$. Of
course, Corollary \ref{BA} says that $\D^1(\mathbb{C}_q(X,Y)
\tensor{\mathbb{C}(X)} \mathbb{C}_q(X,Y), 1 \tensor{\mathbb{C}(X)} 1) = \{ 1 \}$.
Thus, $\D^1$ does not distinguish between the commutative case ($q
= 1$), and the noncommutative case. We propose then the following
definition: given $g \in \gl{\coring{C}}$, we define the
\emph{noncommutative first descent cohomology of $\coring{C}$ at
$g$} as the set of orbits of the action of $U(A)_g$ on
$\gl{\coring{C}}$, notation $N^1(\coring{C},g)$. There is an
obvious surjective map of pointed sets $N^1(\coring{C},g)
\rightarrow \D^1(\coring{C},g)$.
\end{remark}

\begin{example}\label{comoduloalgebra}
Let $H$ be a Hopf algebra over a commutative ring $K$ and consider
any right $H$--comodule algebra $A$ with right coaction $\rho^A :
A \rightarrow A \tensor{K} H$ sending $a \mapsto a_{(0)}\tensor{K}a_{(1)}$ (summation understood). Endow $A \tensor{K} H$ with the
$A$--coring structure given in
\cite[33.2]{Brzezinski/Wisbauer:2003}. Then $1_A \tensor{K} 1_H$
is a grouplike of $A \tensor{K} H$ and
\[
B : = A^{1_A \tensor{K} 1_H} = \{ a \in A : \rho^A(a) = a \tensor{K} 1 \}
\]
Moreover, $1_A \tensor{K} 1_H$ is Galois if and only if $B
\subseteq A$ is a Hopf-Galois $H$--extension. Tomasz Brzezi\'nski
has pointed out \cite[2.6]{Brzezinski:2008} that $\D^1(A \tensor{K}
H, 1_A \tensor{K} 1_H) = \D^1(H,A)$, where the last one refers to
the the first descent cohomology set of $H$ with coefficients in
$A$ defined in \cite{Nuss/Wambst:2007}. We get then the following
consequences of Theorems \ref{simple-cosem} and \ref{mejor}:

\begin{corollary}\label{AH}
Let $B \subseteq A$ be a Hopf-Galois $H$--extension, and assume
that $A$ is left or right IBN. If $B$ is a division ring, then
$\D^1(H,A) = \{ 1 \}$. If, in addition, $B \alpha = \alpha B$ for
every $\alpha \in A$ (e.g., $B \subseteq Center (A))$, then
\begin{equation*}
\gl{A \tensor{K} H} = \left\{\underset{}{} \alpha^{-1} \alpha_{(0)} \tensor{K} \alpha_{(1)} : \alpha \in A^{\times} \right\}
\end{equation*}
is a group with the multiplication
\begin{equation*}
(\alpha^{-1} \alpha_{(0)} \tensor{K} \alpha_{(1)})(\beta^{-1} \beta_{(0)} \tensor{K} \beta_{(1)}) =
\beta^{-1}\alpha^{-1} \alpha_{(0)}\beta_{(0)} \tensor{K} \alpha_{(1)}\beta_{(1)}.
\end{equation*}
We have the isomorphism of groups
\[
\xymatrix{A^{\times}/K^{\times} \ar^-{\cong}[r] & \gl{A \tensor{K} H}, & ( \alpha K^{\times} \ar@{|->}[r] &
\alpha^{-1} \alpha_{(0)} \tensor{K} \alpha_{(1)}) }
\]
\end{corollary}
\end{example}

\begin{example}
A particular case of Example \ref{comoduloalgebra} occurs when $A
= H^a$ the underlying algebra of $H$, and $\rho^H = \Delta$; the comultiplication of $H$. When $K$
is a field, $A^{1 \tensor{K} 1} = K$ and, therefore, $\D^{1}(H,H^a) =
\{ 1 \}$ if $H^a$ is an IBN ring. Moreover, in this case, $\gl{H^a
\tensor{K} H}$  is a group with the multiplication given in
Corollary \ref{AH}. It is easy to check that the map
\begin{equation*}
\xymatrix{
\gl{H} \ar[r] & \gl{H^a \tensor{K} H} & (g \ar@{|->}[r] & 1 \tensor{K} g)
}
\end{equation*}
is a monomorphism of groups. For instance, if $H=K[C_2]$ is the group algebra of the cyclic
group $C_2$ of order $2$ generated by an element $\Sf{g}$.
Then the group of units of the ring $H^a$ is described as follows:
$$U(H^a)\,\,=\,\,\left\{\underset{}{} k + l\Sf{g}|\,\, k,l \in K,
\text{ such that }k^2 - l^2 \,\neq\,0 \right\},$$ the inverse of
$\alpha=k +l\Sf{g} \in U(H^a)$ is given by the element
$\alpha^{-1}=(k^2-l^2)^{-1}(k-l\Sf{g})$. Of course $H^a$ is a
right and left IBN ring, but not a division ring. Therefore,
Corollary \ref{AH} gives a complete description of the group
$\gl{H^a\tensor{K}H}$ which is
$$\left\{\underset{}{} (k^2-l^2)^{-1}\lr{k^2\tensor{K}1 - kl\Sf{g}\tensor{K}1 -l^2\tensor{K}
\Sf{g} + kl\Sf{g}\tensor{K}\Sf{g}}|\, k,l \in K, \text{ such that } k^2-l^2\,\neq\,0\right\}.$$
\end{example}

\providecommand{\bysame}{\leavevmode\hbox
to3em{\hrulefill}\thinspace}
\providecommand{\MR}{\relax\ifhmode\unskip\space\fi MR }
\providecommand{\MRhref}[2]{%
  \href{http://www.ams.org/mathscinet-getitem?mr=#1}{#2}
} \providecommand{\href}[2]{#2}


\begin{thebibliography}{1}

\bibitem{Anderson/Fuller:1974}
F.~W. Anderson and K.~R. Fuller, \emph{Rings and categories of
modules.}, Springer-Verlag, New York, 1974.

\bibitem{Brzezinski:2008}
T.~Brzezi\'{n}ski, \emph{Descent cohomology and corings},
Commun. Algebra  \textbf{36} (2008), 1894-1900.
%
\bibitem{Brzezinski:2002}
T.~Brzezi\'{n}ski, \emph{The structure of corings. {I}nduction
functors, {M}aschke-type theorem, and {F}robenius and
{G}alois-type properties}, Alg.
  Rep. Theory \textbf{5} (2002), 389--410.
%

\bibitem{Brzezinski/Wisbauer:2003}
T.~Brzezi\'{n}ski and R.~Wisbauer, \emph{Corings and comodules},
LMS, vol.
  309, Cambridge University Press, 2003.


%
\bibitem{ElKaoutit/Gomez:2003a}
L.~El~Kaoutit and J.~G\'{o}mez-Torrecillas, \emph{Comatrix
corings: {G}alois corings, {D}escent theory, and a structure
theorem for cosemisimple corings}, Math. Z. \textbf{244} (2003),
887--906.
%

\bibitem{ElKaoutit/Gomez/Lobillo:2004c}
L.~El~Kaoutit, J.~G\'{o}mez-Torrecillas, and F.~J. Lobillo,
\emph{Semisimple  corings}, Alg. Colloq. \textbf{11} (2004), no.~4, 427-442.

\bibitem{Gomez:2002}
J.~G\'{o}mez~Torrecillas, {\em Separable functors in corings},  Int.
J. Math. Math. Sci., {\bf 30}, (2002), 203--225.

\bibitem{Nuss/Wambst:2007}
P. Nuss and M. Wambst, \emph{Non-abelian Hopf cohomology}, J.
Algebra \textbf{312} (2007), 733-754.

\bibitem{Nuss/Wambst:2008}
P. Nuss and M. Wambst, \emph{Non-abelian Hopf cohomology II-The general case}, J.
Algebra \textbf{319} (2008), 4621-4645.
%
\bibitem{Kanzaki:1964}
T.~Kanzaki, \emph{On commutor rings and {G}alois theory of separable algebras},
Osaka J. Math. \textbf{1} (1964), 103--115.
%
\bibitem{Masuoka:1989b}
A.~Masuoka, \emph{Cogalois theory for field extensions}, J. Math.
Soc. Japan \textbf{41} (1989), no.~4, 576--592.
%
%
\bibitem{Serre:1950}
J.~P.~ Serre, Corps locaux. Deuxi\`eme edition, Hermann, Paris,
1968.
%
%
\bibitem{Sweedler:1975}
M.~Sweedler, \emph{The predual theorem to the
{J}acobson-{B}ourbaki theorem}, Trans.
  Amer. Math. Soc. \textbf{213} (1975), 391--406.

\end{thebibliography}
\end{document}